\documentclass{amsart}
\usepackage{amsthm,amsmath,amssymb}
\usepackage{titlesec}
\usepackage{mathrsfs}
\usepackage{mathtools}
\usepackage{comment}
\usepackage{tcolorbox}
\makeatletter
\def\section{\@startsection{section}{1}%
  \z@{.7\linespacing\@plus\linespacing}{.5\linespacing}%
  {\normalfont\scshape\centering}}
\def\subsection{\@startsection{subsection}{2}%
  \z@{.5\linespacing\@plus.7\linespacing}{-.5em}%
  {\normalfont\bfseries}}
\makeatother
\titleformat*{\section}{\large\bfseries}
\titleformat*{\subsection}{\large\bfseries}
\newtheorem{theorem}{Theorem}[section]

\theoremstyle{remark}
\newtheorem{rem}{Remark}
\renewcommand{\Re}{\operatorname{Re}}

\makeatletter
\@namedef{subjclassname@2020}{%
  \textup{2020} Mathematics Subject Classification}
\makeatother
\title[the Riemann zeta function over arithmetic progressions]
{Mean-square values of the Riemann zeta function on arithmetic progressions}
\author{Hirotaka Kobayashi}
\date{}
\address{Graduate School of Mathematics, Nagoya University, Furocho, Chikusaku, Nagoya 464-8602, Japan}
\email{m17011z@math.nagoya-u.ac.jp}
\begin{document}

\begin{abstract}
    We obtain asymptotic formulae for the second discrete moments of the Riemann zeta function over arithmetic progressions $\frac{1}{2} + i(a n + b)$.
    It reveals noticeable relation between the discrete moments and the continuous moment of the Riemann zeta function.
    Especially, when $a$ is a positive integer, main terms of the formula are equal to those for the continuous mean value.
    The proof requires the rational approximation of $e^{\pi k/a}$ for positive integers $k$.
\end{abstract}

\maketitle

\section{Introduction}
In this paper, we shall consider averages
\begin{equation}\label{d-m}
\sum_{an + b \leq T} \left|\zeta \left(\frac{1}{2} + i(a n+b) \right) \right|^2.
\end{equation}
This study is one of the attempts to enrich our knowledge on the vertical distribution of the Riemann zeta function $\zeta(s)$.
Discrete moments of the Riemann zeta function have relation to the distribution of the zeros.
Putnam \cite{Put} showed that there is no infinite arithmetic progression of non-trivial zeros of $\zeta(s)$.
In this direction, there is an important conjecture called the linear independence conjecture, which states that the ordinates of non-trivial zeros of $\zeta(s)$ are linearly independent over $\mathbb{Q}$. In 1942, Ingham \cite{Ing} found out the relation between the linear independence conjecture and the oscillations of $M(x) = \sum_{n \leq x} \mu(n)$, where $\mu(n)$ is the M\"{o}bius function. He showed that the linear independent conjecture implies the failure of the inequality $M(x) \ll x^{1/2}$. For this reason, many mathematicians doubt the inequality.

We may need much more progression to solve the problem. However, we have an easier conjecture in this direction, that is, there are no non-trivial zeros of $\zeta(s)$ in any arithmetic progression. To consider this problem, discrete moments play an important role. Martin and Ng \cite{Ma&Ng} attacked this conjecture for Dirichlet $L$-functions by considering some kinds of discrete means of Dirichlet $L$-functions. Later, Li and Radziwi\l \l \ \cite{Li&Ra} showed that at least one third of the values of the Riemann zeta function on arithmetic progression does not vanish. One of their results (Theorem 2) stated that we have, as $T \to \infty$,
\begin{equation}\label{L-R}
\begin{split}
&\quad
\sum_{n} \left|\zeta \left(\frac{1}{2} + i(an + b) \right) \right|^2 \cdot \phi \left(\frac{n}{T} \right) \\
&=
\int_{\mathbb{R}} \left|\zeta \left(\frac{1}{2} + i(at + b) \right) \right|^2 \cdot \phi \left(\frac{t}{T} \right) dt(1 + \delta(a,b) + o_{a,b,\phi}(1)),
\end{split}
\end{equation}
where $\phi(\cdot)$ is a smooth compactly supported function with support in $[1,2]$, and
\begin{equation*}
  \delta(a,b)
  =
  \begin{dcases*}
    0 & $e^{2\pi k/a}$ is irrational for all $k>0$ \\
    \frac{2\sqrt{rs}\cos(b\log(r/s)) - 2}{rs + 1 - 2\sqrt{rs}\cos(b\log(r/s))} & $e^{2\pi k/a}$ is rational for some $k>0$,
  \end{dcases*}
\end{equation*}
with $r/s \neq 1$ denoting the smallest reduced fraction having a representation in the form $e^{2\pi k/a}$ for some $k>0$.
This clarify the notable correspondence of discrete means to the continuous one.
However, it is difficult to obtain asymptotic formula of the sum (\ref{d-m}), since the error term depends on $\phi$.

\"{O}zbek and Steuding \cite{Oz&St2} proved asymptotic formulae for the first discrete moment of $\zeta(s)$ on certain vertical arithmetic progressions inside the critical strip. The first discrete moment have been studied recently by \"{O}zbek, Steuding and Wegert (see \cite{St&We} and \cite{Oz&St1}). Especially, in \cite{Oz&St1}, they showed that
\begin{equation*}
\lim_{T \to \infty} \frac{1}{T}\sum_{0 \leq n < T} \zeta(s_0 + ina)
=
\begin{dcases}
(1 - l^{-s_0})^{-1} & \text{if} \ a = \frac{2\pi q}{\log l}, \ q \in \mathbb{N}, \ 2 \leq l \in \mathbb{N}, \\
1 & \text{otherwise},
\end{dcases}
\end{equation*}
where $s_0$ may be any complex number with real part in $(0,1)$.

\begin{rem}
  Good \cite{Go} proved asymptotic formulae for fourth moments of the Riemann zeta function on arbitrary arithmetic progressions to the right of the critical line. Namely, he showed that for $\sigma>\frac{1}{2}$
  \begin{equation*}
  \sum_{0 \leq n < T} |\zeta(\sigma + ind)|^4 = T \sum_{m = 1}^{\infty} \frac{d(m)^2}{m^{2\sigma}} + o(T) \quad (T \to \infty),
  \end{equation*}
   where $d$ is not of the form $2\pi l /\log (k_1/k_2)$ with integral $l \neq 0$ and positive integers $k_1 \neq k_2$.
\end{rem}

\begin{rem}
  van Frankenhuijsen \cite{van1} gave an explicit bound for the length of arithmetic progressions of non-trivial zeros of the Riemann zeta function. Later, he \cite{van2} improved the bound.
\end{rem}  

The object of this paper is to prove asymptotic formulae for discrete mean-squares of the Riemann zeta function on vertical arithmetic progressions.
\begin{theorem}\label{main1}
Let $a$ be a real number such that $e^{2\pi k/a}$ is irrational for all positive integer $k$.
We have, as $T \to \infty$,
\begin{equation}\label{th}
\sum_{a n + b \leq T} \left|\zeta \left(\frac{1}{2} + i(a n + b) \right) \right|^2
=
\frac{T}{a}\log T + o_a(T\log T).
\end{equation}
Moreover, when $a$ is a positive integer with $a = o((\log \log T)^{\varepsilon})$, we have
\begin{equation*}
  \sum_{a n + b \leq T} \left|\zeta \left(\frac{1}{2} + i(a n + b) \right) \right|^2
  =
  \frac{T}{a}(\log T + 2\gamma - 1 -\log 2\pi) + O_{A}(a^{-1}T(\log T)^{-A}), 
\end{equation*}
for any fixed $A > 0$.
\end{theorem}
This should be compared with the continuous mean-square
\begin{equation*}
\int_{0}^{T} \left|\zeta \left(\frac{1}{2} + it \right) \right|^2 dt = T \log T + (2\gamma - 1 -\log 2\pi)T + E(T),
\end{equation*}
where $\gamma$ is Euler's constant and $E(T)$ is an error term. Theorem \ref{main1} reveals that the discrete mean values (\ref{th}) equal to the continuous mean value as $T \to \infty$ asymptotically.

Our starting point is the approximate functional equation of $\zeta^2(s)$
\begin{equation}\label{app-zeta2}
\zeta^2(s) = \sum_{n \leq t/2\pi} \frac{d(n)}{n^s} + \chi^2(s) \sum_{n \leq t/2\pi} \frac{d(n)}{n^{1-s}} + R \left(s ; \frac{t}{2\pi} \right),
\end{equation}
where $\chi(s) = 2^s \pi^{s-1} \sin(\pi s/2)\Gamma(1-s)$ and $R(s ; t/2\pi)$ is the error term.
Motohashi \cite{Mo1}, \cite{Mo2} proved that
\begin{equation}\label{R-d-simple}
\chi(1-s) R \left(s ; \frac{t}{2\pi} \right) = -\sqrt{2}\left( \frac{t}{2\pi} \right)^{-1/2}\Delta \left( \frac{t}{2\pi} \right) + O(t^{-1/4}),
\end{equation}
where $\Delta(t/2\pi)$ is the error term in the Dirichlet divisor problem, defined by
\begin{equation*}
\Delta(x) = \sideset{}{'} \sum_{n \leq x} d(n) - x(\log x + 2\gamma - 1) - \frac{1}{4}.
\end{equation*}
Here $\sum'$ indicates that the last term is to be halved if $x$ is an integer. We note that Jutila \cite{Ju} gave another proof of Motohashi's result (\ref{R-d-simple}).

The key step in our proof is to estimate the sum
\begin{equation}\label{d}
  \sum_{1 \le k < (a/2\pi)\log(T/\pi)} e^{\pi k/a} \sum_{m < Te^{-2\pi k/a}}d(m) e(- e^{2\pi k/a} m),
\end{equation}
where $e(x) := \exp(2\pi i x)$.
Bugeaud and Ivi\'{c} \cite{Bug&Iv} also studied a quite similar sum to evaluate the discrete mean value of $E(T)$.
Thus, the same problem arises in the discrete mean value of $E(T)$ and discrete mean-squares of $\zeta(s)$.
They gave the upper bound
\begin{equation*}
  \sum_{1 \le (1/2\pi)\log(T/\pi)} \frac{e^{\pi k}}{k} \sum_{m \leq Te^{-2\pi k}}d(m) e( e^{2\pi k} m)
  \ll T\log T \exp\left(-C\frac{\log \log T}{\log \log \log T}\right),
\end{equation*}
where $C>0$ is some constant.
In our proof, we improve this bound. By the hyperbola method, we obtain the upper bound derived from the exponential sum estimate. The estimate requires a rational approximation of $e^{2\pi k/a}$ by Dirichlet's approximation theorem.
Moreover, in the case of $a$ is a positive integer, we have better estimate applying a result of Waldschmidt \cite{Wal}.
Finally, we apply the argument of Li and Radziwi\l \l \ \cite{Li&Ra} to calculate the sum (\ref{d}).
Consequently, we have
\begin{equation*}
  \begin{split}
    &\quad
    \sum_{1 \le k < (a/2\pi)\log(T/\pi)} e^{\pi k/a} \sum_{m < Te^{-2\pi k/a}}d(m) e(- e^{2\pi k/a} m) \\
    &=
    \begin{dcases}
      o_a(T\log T) & (\text{$a$ is not a integer}), \\
      O_A(T(\log T)^{-A}) & (\text{$a$ is an integer}),
    \end{dcases}
  \end{split}
\end{equation*}
for any fixed $A > 0$.

\begin{rem}
  Bugeaud and Ivi\'{c} \cite{Bug&Iv} have asserted that
  \begin{equation}\label{rem}
    \sum_{n \leq x}E(n) = \pi x + H(x),
  \end{equation}
  where, for some $C>0$, unconditionally
  \begin{equation*}
    H(x) \ll x\log x \exp\left(-C\frac{\log \log x}{\log \log \log x}\right).
  \end{equation*}
  By our improvement of the upper bound of (\ref{d}), this upper bound is also improved to
  \begin{equation*}
    H(x) \ll_{A} x(\log x)^{-A}
  \end{equation*}
  for any fixed $A>0$. Thus we can clarify that the term $\pi x$ in (\ref{rem}) is main term.
  Bugeaud and Ivi\'{c} \cite{Bug&Iv} suggested a conjecture on the upper bound. Now let
  \begin{equation*}
  e^{\pi k} = [a_0(k); a_1(k), \dots]
  \end{equation*}
  be the expansion of $e^{\pi k}$ as a continued fraction for any non-zero integer $k$.
  From the result of Wilton \cite{Wil}, if $a_n(k)$ satisfies $a_n(k) \ll n^{1 + K} \ (K \geq 0)$, then
  \begin{equation*}
  \sum_{m \leq x}d(m) \exp(2\pi i m e^{2\pi k}) \ll x^{1/2} \log^{2 + K} x.
  \end{equation*}
  If this estimate is verified, we can improve the upper bound of $H(x)$ and also Theorem 1.1 with $a = 1$.
\end{rem}

On the other hand, when $e^{2\pi k_0/a}$ is rational for some $k_0$, another main term appears.
\begin{theorem}\label{main2}
  Let $r,s$ be co-prime with $r>2s$.
  Let $a$ be a real number such that 
  \begin{equation*}
    e^{2\pi k_0/a} = \frac{r}{s}
  \end{equation*}
   for some positive integer $k_0$.
We have, as $T \to \infty$,
\begin{equation}\label{th1}
  \begin{split}
  &\quad
  \sum_{a n + b \leq T} \left|\zeta \left(\frac{1}{2} + i(a n+b) \right) \right|^2 \\
  &=
  \frac{T}{a} \log T\left(1 + \frac{2\sqrt{rs}\cos(b\log(r/s))-2}{rs + 1 - 2\sqrt{rs}\cos(b\log(r/s))} + o_{a,b}(1)\right).
  \end{split}
\end{equation}
Moreover, when $k_0 = 1$, we have
\begin{equation}
\begin{split}
  &\quad
  \sum_{a n + b \leq T} \left|\zeta \left(\frac{1}{2} + i(a n+b) \right) \right|^2 \\
  &=
  \frac{T}{a}\left( \log T + 2\gamma - 1 - \log 2\pi \right)\left(1 + \frac{2\sqrt{rs}\cos(b\log(r/s))-2}{rs + 1 - 2\sqrt{rs}\cos(b\log(r/s))} \right) \\
  &\quad
  -\frac{2\sqrt{rs}\cos(b\log(r/s))-2}{rs + 1 - 2\sqrt{rs}\cos(b\log(r/s))}\frac{\sqrt{rs} \log (rs)}{\sqrt{rs} - 1}\frac{T}{a} + o_b(T).
  \end{split}
\end{equation}
\end{theorem}
In this case, the sum (\ref{d}) turns out to be
\begin{equation*}
  \sum_{1 \le k < \log(T/\pi)/\log (r/s)}\sum_{m \le T(s/r)^k} d(m)e\left( - m \left( \frac{r}{s} \right)^k \right).
\end{equation*}
Another main term comes from this sum.

\section{the proof of theorem 1.1}

By (\ref{app-zeta2}), (\ref{R-d-simple}) and the functional equation $\zeta(1-s) = \chi(1-s) \zeta(s)$, we have
\begin{equation*}
\begin{split}
\zeta(s)\zeta(1-s)
&=
\chi(1-s) \sum_{n \leq t/2\pi} \frac{d(n)}{n^s} + \chi (s) \sum_{n \leq t/2\pi} \frac{d(n)}{n^{1-s}} \\
&\quad
-\sqrt{2}\left( \frac{t}{2\pi} \right)^{-\frac{1}{2}}\Delta \left( \frac{t}{2\pi} \right) + O(t^{-1/4}).
 \end{split}
\end{equation*}

It is known that $\Delta(t) \ll t^{1/3+\varepsilon}$. Thus, taking $s = 1/2 + it$, we have
\begin{equation*}
\left|\zeta \left(\frac{1}{2} + it \right) \right|^2 = 2\Re \chi \left(\frac{1}{2} - it \right) \sum_{n \leq t/2\pi} \frac{d(n)}{n^{1/2 + it}} + O(t^{-1/6+\varepsilon}).
\end{equation*}
Hence we consider
\begin{equation*}
\begin{split}
&\quad
\sum_{a n + b \leq T}\left|\zeta \left(\frac{1}{2} + i(a n + b) \right) \right|^2 \\
&=
2 \Re \sum_{a n + b \leq T} \chi \left(\frac{1}{2} - i(a n + b) \right) \sum_{m \leq (a n + b)/2\pi} \frac{d(m)}{m^{1/2 + i(a n + b)}} \\
&\quad
+ O \left(\sum_{a n + b \leq T}(a n + b)^{-1/6+\varepsilon}\right).
\end{split}
\end{equation*}

As for the error term, it is clear that
\begin{equation*}
\sum_{a n + b \leq T} (a n + b)^{-1/6+\varepsilon} \ll a^{-1}T^{5/6}.
\end{equation*}

To consider the main term, we note the following formula
\begin{equation*}
\chi(1-s)=e^{-\pi i/4} \left(\frac{t}{2\pi} \right)^{\sigma - 1/2}\exp \left(it \log \frac{t}{2\pi e} \right)(1+ O (t^{-1}))
\end{equation*}
for fixed $\sigma$ and $t \geq 1$. Using this, we have
\begin{equation*}
    \begin{split}
    &\quad
    \sum_{T < a n + b \leq 2T} \chi \left(\frac{1}{2} - i(a n + b) \right) \sum_{m \leq (a n + b)/2\pi} \frac{d(m)}{m^{1/2 + i(a n + b)}} \\
    &=
    e^{-\pi i/4}\sum_{T < a n + b \leq 2T} \exp \left(i(a n + b) \log \frac{a n + b}{2\pi e} \right)\sum_{m \leq (a n + b)/2\pi} \frac{d(m)}{m^{1/2 + i(a n + b)}} \\
    &\quad
    +O\left( \sum_{T < a n + b \leq 2T} \frac{1}{an+b} \sum_{m\leq (an+b)/2\pi} \frac{d(m)}{m^{1/2}} \right) \\
    &=
    e^{-\pi i/4} \sum_{m \leq T/\pi} \frac{d(m)}{m^{1/2}} \sum_{\max(2\pi m, T) < a n + b \leq 2T} \exp \left(i(a n + b) \log \frac{a n + b}{2\pi em} \right) \\
    &\quad
    + O(a^{-1}T^{1/2} \log T) \\
    &=
    e^{-\pi i/4} \sum_{ m \leq T/2\pi} \frac{d(m)}{m^{1/2}} \sum_{T < a n + b \leq 2T} \exp \left(i(a n + b) \log \frac{a n + b}{2\pi em} \right) \\
    &\quad
    + e^{-\pi i/4} \sum_{T/2\pi < m \leq T/\pi} \frac{d(m)}{m^{1/2}} \sum_{2\pi m < a n + b \leq 2T} \exp \left(i(a n + b) \log \frac{a n + b}{2\pi em} \right) \\
    &\quad
    + O(a^{-1}T^{1/2} \log T) \\
    &=
     e^{-\pi i/4}(S_1 + S_2) + O(a^{-1} T^{1/2} \log T),
    \end{split}
    \end{equation*}
    say. To obtain the second equality, we note that
    \begin{equation*}
    \sum_{m \leq T}\frac{d(m)}{m^{1/2}} \ll T^{1/2}\log T.
    \end{equation*}

    For the convenience, we put
    \begin{equation*}
    f(x) := \frac{a x + b}{2\pi} \log \frac{a x + b}{2\pi em}
    \end{equation*}
    and
    \begin{equation*}
    g_k(x) := f(x) - kx.
    \end{equation*}
    To calculate $S_1$ and $S_2$, we use the saddle-point method.
    
    \subsection{An estimate of $S_1$}
    Since the first and second derivative of $f(x)$ are
    \begin{equation*}
    f'(x) = \frac{a}{2\pi} \log \frac{a x + b}{2\pi m}, \quad f''(x) = \frac{a^2}{2\pi (a x + b)},
    \end{equation*}
    we have (see Proposition 8.7. in \cite{Iw&Ko})
    \begin{equation*}
    \begin{split}
    &\quad
    \sum_{T < a n + b \leq 2T}e^{2\pi if(n)} \\
    &=
    \sum_{(a/2\pi)\log (T/2\pi m) - \theta < k < (a/2\pi)\log (T/\pi m) + \theta} \int_{(T - b)/a}^{(2T - b)/a} e(g_k(x)) dx \\
    &\quad
    + O(\theta^{-1} +\log(a+2)),
    \end{split}
    \end{equation*}
    where $\theta$ is any number with $0 < \theta \leq 1$.
    We choose $\theta = a/2\pi(a + 1)$ and assume $a > 2\pi/\log (T/\pi)$ hereafter.
    If $a \le 2\pi/\log (T/\pi)$, we see that $k = 0$ and then, the integral is $\ll a^{-1}T^{1/2}$ by the second derivative test.

    Here $S_1$ can be rewritten as
    \begin{equation*}
    \begin{split}
    &\quad
    \sum_{m \leq T/2\pi} \frac{d(m)}{m^{1/2}} \sum_{(a/2\pi)\log (T/2\pi m) - \theta < k < (a/2\pi)\log (T/\pi m) + \theta} \int_{(T - b)/a}^{(2T - b)/a} e(g_k(x)) dx \\
    &\quad
    +O(T^{1/2}\log T(\theta^{-1} + \log(a + 2))) \\
    &=
    \sum_{0 \leq k \leq (a/2\pi)\log (T/\pi) + \theta} \\
    &\quad
    \times \sum_{(T/2\pi)e^{-2\pi(k + \theta)/a} < m < (T/\pi)e^{-2\pi(k - \theta)/a}} \frac{d(m)}{m^{1/2}} \int_{(T - b)/a}^{(2T - b)/a} e(g_k(x)) dx \\
    &\quad
    +O(T^{1/2}\log T(\theta^{-1} + \log(a + 2))).
    \end{split}
    \end{equation*}
    We note that the saddle-point of $g_k(x)$, $(2\pi me^{2\pi k/a} - b)/a$, is in $((Te^{-2\pi \theta/a} - b)/a, (2Te^{2\pi \theta/a} - b)/a)$ by the condition in the inner sum.
    We divide the inner sum $((T/2\pi)e^{-2\pi(k + \theta)/a}, (T/\pi)e^{-2\pi(k - \theta)/a})$ into the following five intervals:
    \begin{enumerate}
    \setlength{\itemsep}{0.5cm}
    \item
    $\left(\dfrac{T}{2\pi}e^{-2\pi(k + \theta)/a}, \dfrac{T}{2\pi}e^{-2\pi k/a} - c \right]$,
    \item
    $\left(\dfrac{T}{2\pi}e^{-2\pi k/a} - c, \dfrac{T}{2\pi}e^{-2\pi k/a} + c \right)$,
    \item
    $\left[\dfrac{T}{2\pi}e^{-2\pi k/a} + c, \dfrac{T}{\pi}e^{-2\pi k/a} - c \right]$,
    \item
    $\left(\dfrac{T}{\pi}e^{-2\pi k/a} - c, \dfrac{T}{\pi}e^{-2\pi k/a} + c \right)$,
    \item
    $\left[\dfrac{T}{\pi}e^{-2\pi k/a} + c, \dfrac{T}{\pi}e^{- 2\pi (k - \theta)/a} \right)$,
    \end{enumerate}
    where $c = c(a) := 1/4e(1+a)$. Since $\theta = a/2\pi(a + 1)$, inequalities
    \begin{equation*}
    \frac{T}{2\pi}e^{-2\pi(k + \theta)/a} < \frac{T}{2\pi}e^{-2\pi k/a} - c, \quad \dfrac{T}{2\pi}e^{-2\pi k/a} + c < \dfrac{T}{\pi}e^{-2\pi k/a} - c,
    \end{equation*}
    and
    \begin{equation*}
      \frac{T}{\pi}e^{-2\pi k/a} + c < \frac{T}{\pi}e^{- 2\pi (k - \theta)/a}
    \end{equation*}
    are valid.
    
    (i)
    
    By the first derivative test, we have
    \begin{equation*}
    \left|\int_{(T - b)/a}^{(2T - b)/a}e(g_k(x))dx \right| \leq \cfrac{8\pi}{a \log \cfrac{T}{2\pi me^{2\pi k/a}}},
    \end{equation*}
    and
    \begin{equation*}
    \left|\log \cfrac{T}{2\pi me^{2\pi k/a}} \right| = \left|- \log \left(1 - \frac{T - 2\pi me^{2\pi k/a}}{T}\right) \right| \sim \frac{T - 2\pi me^{2\pi k/a}}{T}.
    \end{equation*}
    Therefore, in this case, the contribution is
    \begin{equation*}
    \begin{split}
    &\ll
    a^{-1}\sum_{0 \leq k \leq (a/2\pi)\log (T/\pi)}\ \sum_{m \leq (T/2\pi)e^{-2\pi k/a} - c} \frac{d(m)}{m^{1/2}} \frac{T}{T - 2\pi me^{2\pi k/a}} \\
    &\ll
    a^{-1}T^{\varepsilon} \sum_{0 \leq k \leq (a/2\pi)\log (T/\pi)}\ \sum_{m \leq (T/2\pi)e^{-2\pi k/a} - c}\frac{T}{2\pi}e^{-2\pi k/a} \cfrac{m^{1/2}}{m\left(\cfrac{T}{2\pi}e^{-2\pi k/a} - m\right)} \\
    &\ll
    a^{-1}T^{\varepsilon} \sum_{0 \leq k \leq (a/2\pi)\log (T/\pi)}\ \sum_{m \leq (T/2\pi)e^{-2\pi k/a} - c} m^{1/2}\left(\frac{1}{m} + \cfrac{1}{\cfrac{T}{2\pi}e^{-2\pi k/a} - m} \right) \\
    &\ll
    a^{-1}T^{1/2 + \varepsilon} \sum_{0 \leq k \leq (a/2\pi)\log (T/\pi)} e^{-\pi k/a} \left(\sum_{m \leq T/2\pi}\frac{1}{m} + a \right) \ll a^{-1}T^{1/2 + \varepsilon}.
    \end{split}
    \end{equation*}

    (v)
    
    In a similar manner, we can see that the contribution of this case is $\ll a^{-1}T^{1/2 + \varepsilon}$.

    (ii), (iv)
    
    First we note that the number of $m$ in each intervals are at most one and
    \begin{equation*}
    m \asymp Te^{-2\pi k/a} \ll T.
    \end{equation*}
    By the second derivative test, we have
    \begin{equation*}
    \left|\int_{(T - b)/a}^{(2T - b)/a}e(g_k(x))dx \right| \leq 16 a^{-1}\sqrt{\pi T}.
    \end{equation*}
    Hence we see that the contribution is
    \begin{equation*}
    \begin{split}
      &\quad
      a^{-1}T^{1/2} \sum_{0 \leq k \leq (a/2\pi)\log (T/\pi) + \theta}\ \sum_{\substack{(T/2\pi)e^{-2\pi k/a} - c < m < (T/2\pi)e^{-2\pi k/a} + c \\ (T/\pi)e^{-2\pi k/a} - c < m < (T/\pi)e^{-2\pi k/a} + c}} \frac{d(m)}{m^{1/2}} \\
      &\ll
      e^{1/2a}T^{1/2+\varepsilon} \ll T^{5/6}.
    \end{split}
    \end{equation*}

    (iii)
    
    In this case, the sum of $k$ starts from $1$. If $k = 0$, then we have that $T/2\pi + c < m$. However, this is impossible, because we consider the case $m \leq T/2\pi$.
    Using the saddle-point method (see for example Corollary 8.15. in \cite{Iw&Ko}), we have
    \begin{equation*}
    \begin{split}
    \int_{(T - b)/a}^{(2T - b)/a}e(g_k(x))dx 
    &=
    e^{\pi i/4}\frac{2\pi}{a} e^{\pi k/a + 2\pi i bk/a} \sqrt{m}\exp\left(-2\pi mi e^{2\pi k/a} \right) \\
    &\quad
    +O \left(\frac{T}{2T - 2\pi me^{2\pi k/a}} + \frac{T}{2\pi me^{2\pi k/a} - T} + 1 \right).
    \end{split}
    \end{equation*}
    By the same argument with the case (i) and (v), we see that the contribution of the error term is $\ll a^{-1}T^{1/2 + \varepsilon}$.
    Finally, we have to consider
    \begin{equation*}
      \begin{split}
        &\quad
        \sum_{1 \leq k < (a/2\pi)\log (T/\pi) + \theta} e^{(\pi + 2\pi i b)k/a} \\
        &\quad
        \times \sum_{(T/2\pi)e^{-2\pi k/a} + c \le m \le (T/\pi)e^{-2\pi k/a} - c} d(m) e(-e^{2\pi k/a}m).
      \end{split}
    \end{equation*}
    However we see that
    \begin{equation*}
      \sum_{1 \leq k < (a/2\pi)\log (T/\pi) + \theta} e^{(\pi + 2\pi i b)k/a} \sum_{\substack{(T/2\pi)e^{-2\pi k/a} + c \le m \le (T/2\pi)e^{-2\pi k/a} \\ (T/\pi)e^{-2\pi k/a} \le m \le (T/\pi)e^{-2\pi k/a} - c}} \ll T^{1/2 + \varepsilon},
    \end{equation*}
    and when $(a/2\pi)\log (T/\pi) \le k < (a/2\pi)\log (T/\pi) +\theta$, the inner sum is empty sum. Thus we calculate
    \begin{equation}\label{irrational-sum}
      \sum_{1 \leq k < (a/2\pi)\log (T/\pi)} e^{(\pi + 2\pi i b)k/a} \sum_{(T/2\pi)e^{-2\pi k/a} < m < (T/\pi)e^{-2\pi k/a}} d(m) e(-e^{2\pi k/a}m).
    \end{equation}

    \textit{\textbf{Case 1}}
    
    Suppose that $e^{2\pi k/a}$ is irrational for all positive integer $k$.
    In this case, we consider the sum
    \begin{equation*}
    \sum_{m \leq M} d(m)e(\alpha m),
    \end{equation*}
    where $\alpha = \alpha (k) = e^{2 \pi k/a}$. By the hyperbola method, we have
    \begin{equation*}
    \begin{split}
    \sum_{m \leq M} d(m)e(\alpha m)
    &=
    \sum_{uv \leq M} e(\alpha uv) \\
    &=
    2\sum_{u \leq \sqrt{M}}\sum_{u < v <M/u}e(\alpha uv) + \sum_{u \leq \sqrt{M}} e(\alpha u^2).
    \end{split}
    \end{equation*}
    We can plainly see that the second term is $\ll \sqrt{M}$.
    
    As for the first sum, making use of the estimate
    \begin{equation*}
    \sum_{N_1 < n \leq N_2}e(\alpha n) \ll \min (N_2 - N_1, ||\alpha||^{-1}),
    \end{equation*}
    we obtain
    \begin{equation*}
    \sum_{u \leq \sqrt{M}}\sum_{u < v <M/u}e(\alpha uv) \ll \sum_{u \leq \sqrt{M}}\min \left( \frac{M}{u}, ||\alpha u||^{-1}\right).
    \end{equation*}
    Here, by Dirichlet's approximation theorem, we can take integers $p = p(k), q = q(k)$ such that $(p,q) = 1$, $1 \leq q \leq \sqrt{M}$ and
    \begin{equation}\label{approximation}
    \left|\alpha - \frac{p}{q} \right| \leq \frac{1}{q\sqrt{M}} \leq \frac{1}{q^2}.
    \end{equation}
    In this situation, we have
    \begin{equation*}
    \sum_{u \leq \sqrt{M}}\min \left( \frac{M}{u}, ||\alpha u||^{-1}\right) \ll M\left(\frac{1}{q} + \frac{1}{\sqrt{M}} + \frac{q}{M} \right)\log (qM).
    \end{equation*}
    Applying this estimate with $M = T e^{-2\pi k/a}$, we see that the sum (\ref{irrational-sum}) is
    \begin{equation*}
    \begin{split}
    &
    \ll \sum_{1 \leq k < (a/2\pi)\log (T/2\pi)} (T e^{-\pi k/a} q(k)^{-1} + T^{1/2} + e^{\pi k/a}q(k) )\log T \\
    &\ll
    T\log T\sum_{1 \leq k < (a/2\pi)\log (T/2\pi)} e^{-\pi k/a} q(k)^{-1} + a T^{1/2} \log^2 T.
    \end{split}
    \end{equation*}
    We note that the condition (\ref{approximation}) leads to $p \asymp q\alpha$ and $pq  \gg \alpha$. Therefore we can find an $A$ such that
    \begin{equation*}
    \begin{split}
    \sum_{1 \leq k < (a/2\pi)\log (T/2\pi)} e^{-\pi k/a} q(k)^{-1}
    &\asymp
    \sum_{1 \leq k < (a/2\pi)\log (T/2\pi)} (p(k)q(k))^{-1/2} \\
    &\ll
    \sum_{1 \leq k < A} (p(k)q(k))^{-1/2} + \sum_{A \leq k} e^{-\pi k/a} \\
    &\ll
    \sum_{1 \leq k <A} (p(k)q(k))^{-1/2} + e^{-C\log \log T}.
    \end{split}
    \end{equation*}
    By the condition (\ref{approximation}), when $T$ tends to infinity, for each $k$, $p(k)q(k)$ does so. Thus
    \begin{equation*}
        \sum_{1 \leq k < A} (p(k)q(k))^{-1/2} = o_{a}(1)
    \end{equation*}
    as $T \to \infty$.
    We conclude that the sum (\ref{irrational-sum}) is $o_{a}(T\log T)$ and so as $S_1$.
    
\textit{\textbf{Case 2}}

    Moreover, if $a$ is a positive integer, then we can obtain better estimate, applying an inequality
    \begin{equation*}
    \left|e^{\pi k/a} - \frac{p}{q} \right| > \exp \{ -2^{72}(\log 2k)(\log 2a)(\log p)(\log \log p) \}
    \end{equation*}
    due to Waldschmidt. By this bound, when the condition (\ref{approximation}) is valid, we have
    \begin{equation*}
    e^{\pi k/a}T^{-1/2} \geq \exp \{-c(\log 2k)(\log 2a)(\log p)(\log \log p)\}.
    \end{equation*}
    When $k < (\log \log T)^{1+\varepsilon}$, we obtain $(\log p)(\log \log p) \gg \log T/(\log 2a)(\log \log T)^{\varepsilon}$,
    by the above inequality, and hence $\log p \gg \log T/ (\log 2a)(\log \log T)^{1+\varepsilon}$. Therefore we find that
    \begin{equation*}
    \begin{split}
    \sum_{1 \leq k <(\log \log T)^{1+\varepsilon}} (p(k)q(k))^{-1/2}
    &\asymp 
    \sum_{1 \leq k <(\log \log T)^{1+\varepsilon}} (p(k)q(k))^{-1/4} e^{\pi k/2a} e^{-(\log a)/2} \\
    &\ll
    e^{-c\log T/(\log 2a)(\log \log T)^{1+\varepsilon}}  \sum_{1 \leq k <(\log \log T)^{1+\varepsilon}} 1 \\
    &\ll
    e^{-c\log T/(\log 2a)(\log \log T)^{1+\varepsilon}} \\
    &\ll
    e^{-c\log T/(\log \log T)^{1+\varepsilon}}.
    \end{split}
    \end{equation*}
    Thus, in this case, $S_1$ is $O(a^{-1}T(\log T)^{-A})$ for any fixed $A>0$.
        
    \subsection{A contribution of $S_2$}
    As in the case $S_1$, we have
    \begin{equation*}
      \begin{split}
    &\quad    
    \sum_{2\pi m < a n + b \leq 2T}e^{2\pi if(n)} \\
    &=
    \sum_{0 \leq k < (a/2\pi)\log (T/\pi m) + \theta} \int_{(2\pi m-b)/a}^{(2T-b)/a} e(g_k(x)) dx + O(\theta^{-1} + \log (a\log T)),
      \end{split}
    \end{equation*}
    where $\theta = a/2\pi(a + 1)$.
    
    Since
    \begin{equation*}
    \sum_{m \leq T} \frac{d(m)}{m^{1/2}} \ll T^{1/2}\log T,
    \end{equation*}
    the contribution of the error term is $\ll T^{1/2} \log T \log (a\log T)$.
    
    Finally, we calculate
    \begin{equation*}
    \sum_{T/2\pi < m \leq T/\pi} \frac{d(m)}{m^{1/2}} \sum_{0 \leq k < (a/2\pi)\log (T/\pi m) + \theta} \int_{(2\pi m-b)/a}^{(2T-b)/a} e(g_k(x)) dx.
    \end{equation*}
    When $k = 0$, the above is
    \begin{equation*}
    \sum_{T/2\pi < m \leq T/\pi} \frac{d(m)}{m^{1/2}} \int_{(2\pi m-b)/a}^{(2T-b)/a} e(f(x)) dx.
    \end{equation*}

    The saddle-point of $f(x)$ is $(2\pi m - b)/a$, thus, the saddle-point method leads to
    \begin{equation*}
    \int_{(2\pi m-b)/a}^{(2T-b)/a} e(f(x)) dx = e^{\pi i/4} \frac{\pi}{a} \sqrt{m} + O \left(\frac{T}{2T - 2\pi m} + 1 \right).
    \end{equation*}
    Therefore 
    \begin{equation}\label{main}
    \begin{split}
      &\quad
      \sum_{T/2\pi < m \leq T/\pi} \frac{d(m)}{m^{1/2}} \int_{(2\pi m-b)/a}^{(2T-b)/a} e(f(x)) dx \\
      &= e^{\pi i/4} \frac{\pi}{a} \sum_{T/2\pi < m \leq T/\pi} d(m) + O(T^{1/2 + \varepsilon}).
    \end{split}
    \end{equation}

    When $k \neq 0$, we have to calculate
    \begin{equation*}
    \begin{split}
    &\quad
    \sum_{T/2\pi < m \leq T/\pi} \frac{d(m)}{m^{1/2}} \sum_{1 \leq k < (a/2\pi)\log (T/\pi m) + \theta} \int_{(2\pi m-b)/a}^{(2T-b)/a} e(g_k(x)) dx \\
    &=
    \sum_{1 \leq k < (a/2\pi)\log 2} \ \sum_{T/2\pi < m \leq (T/\pi)e^{-2\pi(k - \theta)/a}} \frac{d(m)}{m^{1/2}} \int_{(2\pi m-b)/a}^{(2T-b)/a} e(g_k(x)) dx,
    \end{split}
    \end{equation*}
    but in a similar manner to $S_1$, we can see that the contribution of this case is
    \begin{equation*}
      \begin{cases}
        o_a(T\log T) & \text{($a$ is not integer)}, \\
        O_A(a^{-1}T(\log T)^{-A}) & \text{($a$ is integer)}.
      \end{cases}
    \end{equation*}     

    \subsection{Conclusion}
    From the above, we can see that
    \begin{equation*}
    \begin{split}
    &\quad
    \sum_{T < a n + b \leq 2T} \chi \left(\frac{1}{2} - i(a n + b) \right) \sum_{m \leq (a n + b)/2\pi} \frac{d(m)}{m^{1/2 + i(a n + b)}} \\
    &=
    \frac{\pi}{a} \sum_{T/2\pi < m \leq T/\pi} d(m) + R(T),
    \end{split}
    \end{equation*}
    and so,
    \begin{equation*}
    \sum_{T < a n + b \leq 2T}\left|\zeta \left(\frac{1}{2} + i(a n + b) \right) \right|^2
    =
    \frac{2\pi}{a} \sum_{T/2\pi < m \leq T/\pi} d(m) + R(T),
    \end{equation*}
    where
    \begin{equation*}
        R(T)
        =
        \begin{cases}
          o_{a}(T\log T) & \text{($a$ is not integer)} \\
          O_A(a^{-1}T(\log T)^{-A}) & \text{($a$ is integer)}.
        \end{cases}
    \end{equation*}
    
    Finally, replacing $T$ by $T/2, T/4$, and so on, and adding we have
    \begin{equation*}
    \begin{split}
    \sum_{an+b \leq T}\left|\zeta \left(\frac{1}{2} + i(a n + b) \right) \right|^2
    &=
    \frac{2\pi}{a} \sum_{ m \leq T/2\pi} d(m) + R(T) \\
    &=
    \frac{T}{a}(\log T +2\gamma-1-\log 2\pi) + R(T).
    \end{split}
    \end{equation*}

    This completes the proof.    
    
\section{the proof of theorem 1.2}
Now we consider the case that $e^{2\pi k_0/a}$ is rational for some $k_0$. In this case, we can write
\begin{equation*}
a = \frac{2\pi k_0}{\log (r/s)}
\end{equation*}
with relatively prime integers $r$ and $s$ and $|r|$ minimal. Let $l$ be the maximal positive integer such that $r/s=(x/y)^l$ with $r$, $s$ relatively prime.
Then,
\begin{equation*}
a = \frac{k_0}{l} \frac{2\pi}{\log (x/y)}.
\end{equation*}

For each $k$ divisible by $k_0$, $e^{2\pi k/a} = (r/s)^{k/k_0}$ is rational. On the other hand, $k$ is not divisible by $k_0$, $e^{2\pi k/a} = (x/y)^{kl/k_0}$ is irrational since $k_0 \nmid l$.
Therefore, the sum (\ref{irrational-sum}) can be divided as
\begin{equation*}
\begin{split}
&\quad    
\sum_{\substack{1 \leq k < (a/2\pi)\log (T/\pi) \\ k_0 \mid k}} e^{(\pi + 2\pi i b)k/a} \sum_{(T/2\pi)e^{-2\pi k/a} < m < (T/\pi)e^{-2\pi k/a}} d(m) e(-e^{2\pi k/a}m) \\
&
+ \sum_{\substack{1 \leq k < (a/2\pi)\log (T/\pi) \\ k_0 \nmid k}} e^{(\pi + 2\pi i b)k/a} \sum_{(T/2\pi)e^{-2\pi k/a} < m < (T/\pi)e^{-2\pi k/a}} d(m) e(-e^{2\pi k/a}m).
\end{split}
\end{equation*}
When $k_0=1$, the second sum is empty sum, otherwise it is $o_{a}(T\log T)$ as can be seen by repeating the same argument as in the Case 1 of the proof of Theorem \ref{main1}.

As for the first sum, we separate the outer sum into two parts as following
\begin{equation*}
\begin{split}
&\quad
\sum_{1 \leq k < (a/2\pi k_0)\log (T/\pi)} \left(\frac{r}{s} \right)^{(1/2+ib)k} \sum_{(T/2\pi)(s/r)^{k} < m < (T/\pi)(s/r)^{k}}d(m) e \left(-m \left(\frac{r}{s} \right)^{k} \right) \\
&=
\sum_{1 \leq k < \log (T/\pi)/\log (rs)} \left(\frac{r}{s} \right)^{(1/2+ib)k} \sum_{(T/2\pi)(s/r)^{k} < m < (T/\pi)(s/r)^{k}}d(m) e \left(-m \left(\frac{r}{s} \right)^{k} \right) \\
&\quad
+ \sum_{ \log (T/\pi)/\log (rs) \leq k < \log (T/\pi)/\log (r/s)} \left(\frac{r}{s} \right)^{(1/2+ib)k} \\
&\quad
\times \sum_{(T/2\pi)(s/r)^{k} < m < (T/\pi)(s/r)^{k}}d(m) e \left(-m \left(\frac{r}{s} \right)^{k} \right).
\end{split}
\end{equation*}
Using
\begin{equation*}
  \sum_{m \le x} d(m)e \left( \frac{mr}{s} \right) = \frac{x}{s}\left( \log x + 2\gamma - 1 - 2\log s \right) + O((\sqrt{x} + s)\log 2s),
\end{equation*}
we have
\begin{equation*}
\begin{split}
  &\quad
\sum_{1 \leq k < \log (T/\pi)/\log (rs)} \left(\frac{(r/s)^{ib}}{\sqrt{rs}}\right)^{k}\frac{T}{2\pi} \left(\log \frac{T}{2\pi} + 2\gamma - 1 -k\log (rs) \right) \\
&\quad
+O_b(T^{1/2} \log^2 T) \\
&=
\frac{(r/s)^{ib}}{\sqrt{rs} - (r/s)^{ib}}\frac{T}{2 \pi}\left( \log \frac{T}{2\pi} + 2\gamma - 1 - \frac{\sqrt{rs} \log (rs)}{\sqrt{rs} - 1} \right) +O_b(T^{1/2} \log^2 T).
\end{split}
\end{equation*}
As for another sum, since
\begin{equation*}
  \sum_{m \le x}d(m) \ll x\log x,
\end{equation*}
we see that
\begin{equation*}
\begin{split}
  &\quad
  \sum_{ \log (T/\pi)/\log (rs) \leq k < \log (T/\pi)/\log (r/s)} \left(\frac{r}{s} \right)^{(1/2+ib)k} \\
  &\quad
  \times \sum_{(T/2\pi)(s/r)^{k} < m < (T/\pi)(s/r)^{k}}d(m) e \left(-m \left(\frac{r}{s} \right)^{k} \right) \\
  &\ll
  T\log T \sum_{ \log (T/\pi)/\log (rs) \leq k < \log (T/\pi)/\log (r/s)} \left( \frac{r}{s} \right)^{-k/2} \ll T^{1- \log(r/s)/2\log(rs)}\log T.
\end{split}  
\end{equation*}
Therefore, we obtain
\begin{equation*}
\begin{split}
  &\quad
  S_1 \\
  &=
  e^{\pi i/4}\frac{(r/s)^{ib}}{\sqrt{rs} - (r/s)^{ib}}
  \begin{dcases}
    \frac{T}{a} \log T(1+o_{a,b}(1)) & (k_0 > 1) \\
    \frac{T}{a}\left( \log \frac{T}{2\pi} + 2\gamma - 1 - \frac{\sqrt{rs} \log (rs)}{\sqrt{rs} - 1} \right) + o_{b}(T) & (k_0=1).
  \end{dcases}
\end{split}
\end{equation*}

As for $S_2$, we have
\begin{equation*}
  \begin{split}
    S_2
    &=
    e^{\pi i/4} \frac{\pi}{a} \sum_{T/2\pi < m \leq T/\pi} d(m) + O(T^{1/2 + \varepsilon}) \\
    &\quad
    +\sum_{1 \leq k < (a/2\pi)\log 2} \ \sum_{T/2\pi < m \leq (T/\pi)e^{-2\pi(k - \theta)/a}} \frac{d(m)}{m^{1/2}} \int_{(2\pi m-b)/a}^{(2T-b)/a} e(g_k(x)) dx.
  \end{split}
\end{equation*}
When $k_0=1$, by the condition of $r$ and $s$, the double sum is empty, otherwise it is $o_{a}(T\log T)$ as can be seen by repeating the same argument as in the Case 1 of the proof of Theorem \ref{main1}.

\section*{Acknowledgement}

The author would like to thank Professors Kohji Matsumoto and Yoonbok Lee
for their helpful advice. The author is supported by Foundation of Research Fellows,
The Mathematical Society of Japan.

\end{document}